\documentclass{amsart}

\newtheorem{theorem}{Theorem}[section]

\newtheorem{corollary}{Corollary}[section]
\newtheorem{lemma}{Lemma}[section]
\newtheorem{note}{Note}[section]

\usepackage{graphicx}
\usepackage{subfigure}
\usepackage{amsfonts}
\usepackage{amssymb}
\usepackage{amsmath,latexsym,amsthm}

\theoremstyle{definition}

\theoremstyle{remark}

\begin{document}

\title[The Glasser-Manna-Oloa integral]{The Laplace transform of the digamma 
function: an integral due to Glasser, Manna and Oloa}

\author{Tewodros Amdeberhan}
\address{Department of Mathematics,
Tulane University, New Orleans, LA 70118}
\email{tamdeber@tulane.edu}

\author{Victor H. Moll}
\address{Department of Mathematics,
Tulane University, New Orleans, LA 70118}
\email{vhm@math.tulane.edu}

\thanks{The work of the second author was partially funded by
$\text{NSF-DMS } 0409968$.}



\numberwithin{equation}{section}

\newcommand{\imagpart}{\mathop{\rm Im}\nolimits}
\newcommand{\realpart}{\mathop{\rm Re}\nolimits}
\newcommand{\no}{\noindent}

\begin{abstract}
The definite integral
\begin{equation}
M(a):= \frac{4}{\pi} \int_{0}^{\pi/2} \frac{x^{2} \, dx }
{x^{2} + \ln^{2}( 2 e^{-a} \cos x ) },
\nonumber
\end{equation}
\noindent
is related to the Laplace transform of the digamma function
\begin{equation}
L(a) := \int_{0}^{\infty} e^{-a s} \psi(s+1) \, ds,
\nonumber
\end{equation}
\noindent
by $M(a) = L(a) + \gamma/a$ when $a > \ln 2$. We establish an analytic 
expression for $M(a)$ in the complementary range $0 < a \leq \ln 2$. 
\end{abstract}

\maketitle

\section{Introduction} \label{sec-intro}
\setcounter{equation}{0}

The classical table of integrals by I. S. Gradshteyn and I. M. Ryzhik \cite{gr}
contains a large collection  organized in sections 
according to the form of the integrand. In each section  one finds
significant variation on the complexity of the integrals. For example, 
section $4.33-4.34$,
with the title {\em Combinations of logarithms and exponentials}, presents
the elementary formula $4.331.1$: for $ a > 0$,
\begin{equation}
\int_{0}^{\infty} e^{-ax} \ln x \, dx = -\frac{\gamma + \ln a }{a},
\label{int-00}
\end{equation}
\noindent
where $\gamma$ is the {\em Euler constant}
\begin{equation}
\gamma = \lim\limits_{n \to \infty} \sum_{k=1}^{n} \frac{1}{k} - \ln n,
\end{equation}
\noindent
as well as the more elaborate $4.332.1$ and $4.325.6$:
\begin{equation}
\int_{0}^{\infty} \frac{\ln x \, dx}{e^{x}+e^{-x} -1} = 
\int_{0}^{1} \ln \ln \left( \frac{1}{x} \right) \frac{dx}{x^{2}-x+1} = 
\frac{2 \pi}{\sqrt{3}} \left( \frac{5}{6} \ln 2 \pi - \ln \Gamma \left( 
\frac{1}{6} \right) \right). 
\nonumber
\end{equation}
\noindent
The difficult involved in the evaluation of a definite integral is hard to 
measure from the complexity of the integrand. For instance, the evaluation 
of {\em Vardi's integral}, 
\begin{equation}
\int_{\pi/4}^{\pi/2} \ln \ln \tan x \, dx = 
\int_{0}^{1} \ln \ln \left( \frac{1}{x} \right) \frac{dx}{1+x^{2}} = 
\frac{\pi}{2} 
\ln \left( \frac{\Gamma( \tfrac{3}{4}) \, \sqrt{2 \pi}}{\Gamma( \tfrac{1}{4}) }
\right),
\end{equation}
\noindent
that appears as $4.229.7$ in \cite{gr}, requires a reasonable amount 
of Number Theory. The second form is $4.325.4$, found in the section 
entitled {\em Combinations of logarithmic functions of more
complicated arguments and powers}. 
The reader will find in \cite{vardi1} a discussion of this 
formula.

It is a remarkable 
fact that combinations of elementary functions in the integrand often 
exhibits 
definite integrals whose evaluations are far from elementary. We have initiated 
a systematic study of the formulas in \cite{gr} in 
the series \cite{moll-gr5, moll-gr9, moll-gr1, moll-gr2, 
moll-gr3, moll-gr4}. The papers are organized according to the 
{\em combinations} appearing in the integrand. Even the elementary cases, such 
as the combination of logarithms and rational function discussed in 
\cite{moll-gr9} entail interesting results. The evaluations
\begin{eqnarray}
\int_{0}^{b} \frac{\ln t \, dt}{(1+t)^{n+1}} & =  &  
\frac{1}{n} \left[ 1 - (1+b)^{-n} \right] \ln b - \frac{1}{n}  \ln(1+b)
\label{stir1}  \\
& - & \frac{1}{n (1+b)^{n-1}} \sum_{j=1}^{n-1} 
\frac{1}{j!} \binom{n-1}{j} |s(j+1,2)| b^{j}, \nonumber
\end{eqnarray}
\noindent
for $b > 0$ and $n \in \mathbb{N}$ produces an explicit formula for the 
case where the rational function has a single pole. Here $s(n,k)$ are the 
{\em Stirling numbers of the first kind} counting the number of  
permutations of 
$n$ letters having exactly $k$ cycles. 
The case of a purely imaginary pole is expressed in terms of the 
rational function
\begin{equation}
p_{n}(x) = 
\sum_{j=1}^{n} \frac{2^{2j}}{2j \binom{2j}{j}} \frac{x}{(1+x^{2})^{j}},
\end{equation}
\noindent
as
\begin{equation}
\int_{0}^{x} \frac{\ln t \, dt}{(1+t^{2})^{n+1}}  
= \frac{\binom{2n}{n}}{2^{2n}} \left[ g_{0}(x) + p_{n}(x) \ln x - 
\sum_{k=0}^{n-1} 
\frac{\tan^{-1} x  + p_{k}(x)}{2k+1} \right],
\nonumber
\end{equation}
\noindent
with
\begin{equation}
g_{0}(x) = \ln x \, \tan^{-1}x - \int_{0}^{x} \frac{\tan^{-1}t}{t} \, dt.
\end{equation}
\noindent
The special case $x=1$ becomes
\begin{equation}
\int_{0}^{1} \frac{\ln t \, dt}{(1+t^{2})^{n+1}} = -2^{-2n} \binom{2n}{n} 
\left( G + \sum_{k=0}^{n-1} \frac{ \tfrac{\pi}{4} + p_{k}(1)}{2k+1} \right),
\end{equation}
\noindent 
where 
\begin{equation}
G = \sum_{k=0}^{\infty} \frac{(-1)^{k}}{(2k+1)^{2}}
\end{equation}
\noindent
is the {\em Catalan's constant}. The values 
\begin{equation}
p_{k}(1) =  \sum_{j=1}^{k} \frac{2^{j}}{2j \binom{2j}{j} }
\end{equation}
\noindent
do not admit a closed-form (in the sense of \cite{aequalsb}), but 
they do satisfy the three term recurrence 
\begin{equation}
(2k+1)p_{k+1}(1) - (3k+1)p_{k}(1) + k p_{k-1}(1) = 0. 
\end{equation}

%

The study of definite integrals where the integrand is a combination of 
powers, logarithms and trigonometric functions was initiated by Euler 
\cite{euler72}, with the evaluation of 
\begin{equation}
\int_{0}^{\pi/2} x \, \ln (2 \cos x ) \, dx  = 
- \frac{7}{16} \zeta(3),
\end{equation}
\noindent
and
\begin{equation}
\int_{0}^{\pi/2} x^{2} \, \ln (2 \cos x ) \, dx = 
- \frac{\pi}{4} \zeta(3),
\end{equation}
\noindent
that appear in his study of the {\em Riemann zeta} function at the odd 
integers. 
These type of integrals have been investigated in \cite{kolbig2}, 
\cite{yue2}.
The {\em intriguing  integral} \cite{borwein95a},
\begin{equation}
\int_{0}^{\pi/2} x^{2} \, \ln ^{2}(2 \cos x ) \, dx = 
\frac{11 \pi}{16} \zeta(4) = \frac{11 \pi^{5}}{1440},
\end{equation}
\noindent
was first conjectured on the basis of a numerical computation 
by Enrico Au-Yueng, while an undergraduate student at the University of 
Waterloo.  \\

Recently O. Oloa \cite{oloa1} considered the integral
\begin{equation}
M(a):= \frac{4}{\pi} \int_{0}^{\pi/2} \frac{x^{2} \, dx }
{x^{2} + \ln^{2}( 2 e^{-a} \cos x ) },
\label{oloa-11}
\end{equation}
\noindent
and later established the value 
\begin{equation}
M(0)= \frac{4}{\pi} \int_{0}^{\pi/2} \frac{x^{2} \, dx }
{x^{2} + \ln^{2}( 2  \cos x ) } = \frac{1}{2}(1 + \ln( 2 \pi) - \gamma).
\label{value-m0}
\end{equation}

Oloa's method of proof relies on the expansion
\begin{equation}
\frac{x^{2}}{x^{2} + \ln^{2}( 2 \cos x)} = x \sin 2x + 
\sum_{n=1}^{\infty} (-1)^{n-1} \left( \frac{a_{n}}{n!} - 
\frac{a_{n+1}}{(n+1)!} \right) x \sin( 2 n x),
\end{equation}
\noindent
where 
\begin{equation}
a_{n} := \int_{0}^{1} (t)_{n} \, dt, 
\end{equation}
\noindent
with $(t)_{n} = t(t+1) \cdots (t+n-1)$ the {\em Pochhammer symbol}. The 
standard relation
\begin{equation}
(t)_{n} = \sum_{k=1}^{n} |s(n,k)| t^{k}, 
\end{equation}
\noindent
gives 
\begin{equation}
a_{n} =  \sum_{k=1}^{n} \frac{|s(n,k)|}{k+1}.
\end{equation}

M. L. Glasser and D. Manna \cite{glasser-manna1} introduced 
the function 
\begin{equation}
L(a) := \int_{0}^{\infty} e^{-a s} \psi(s+1) \, ds,
\label{glasser-2}
\end{equation}
\noindent
where $\psi(x) = \frac{d}{dx} \ln \Gamma(x)$ is the {\em digamma 
function}.
Integration by parts and using (\ref{int-00}) gives 
\begin{equation}
L(a) = -\gamma - \ln a + a \int_{0}^{\infty} e^{-at} \ln \Gamma(t) \, dt.
\end{equation}
\noindent
The main result in \cite{glasser-manna1} gives a relation between $M(a)$
and $L(a)$.  

\begin{theorem}
\label{gm-1}
Assume that $a > \ln 2$. Then
\begin{equation}
M(a) = L(a) + \frac{\gamma}{a}.
\label{gm-11}
\nonumber
\end{equation}
\noindent
That is, for $a > \ln2 $, 
\begin{equation}
M(a)  = \frac{\gamma}{a} - \gamma - \ln a + 
a \int_{0}^{\infty} e^{-at} \ln \Gamma(t) \, dt.
\end{equation}
\end{theorem}

The proof in \cite{glasser-manna1} begins with the representation 
\begin{equation}
\int_{0}^{\pi/2} \cos^{\nu}x \, \cos a x \, dx = 
\frac{\pi \Gamma(\nu+2)}{2^{\nu+1} (\nu+1) \, \Gamma( 1 + \tfrac{\nu}{2} +
\tfrac{a}{2} ) \, \Gamma( 1 + \tfrac{\nu}{2} - \tfrac{a}{2} )},
\end{equation}
\noindent
borrowed from $3.621.9$ in \cite{gr}. Differentiating with respect to $a$, 
evaluating at $a = s$, and using $\psi(1) = - \gamma$ yields 
\begin{equation}
\psi(s+1) = \frac{2^{s+2}}{\pi} 
\int_{0}^{\pi/2} x \cos^{s}x \, \sin(sx) \, dx - \gamma.
\label{psi-1}
\end{equation}
\noindent
Replacing in (\ref{glasser-2}) produces
\begin{equation}
L(a) + \frac{\gamma}{a} = 
-\frac{4}{\pi} \imagpart{\int_{0}^{\infty} \int_{0}^{\pi/2} 
x e^{s( \ln \left[2 e^{-a} \cos x \right] - i x)} \, dx \, ds }.
\end{equation}
\noindent
The identity (\ref{gm-11}) now is an immediate consequence of the 
$s$-integral:
\begin{equation}
\int_{0}^{\infty} e^{s( \ln \left[2 e^{-a} \cos x \right] - i x} \, ds
= \frac{1}{ix - \ln \left[ 2 e^{-a} \cos x \right]}.
\end{equation}

The authors also succeed in a series expansion of $M(a)$ while 
they recognize as
a hypergeometric function in two variables, and noted (quoting from 
\cite{glasser-manna1}) {\em strongly suggests that for general value of $a$, 
no further progress is possible}.  This hypergeometric interpretation led
the authors \cite{glasser-manna1} to
\begin{equation}
M(0) = 1+ \frac{1}{2} \int_{0}^{1} t(1-t) \,{_{3}}F_{2}(1,1,2-t; 2,3; 1) \, dt 
\end{equation}
\noindent 
for which they invoke
\begin{equation}
{_{3}}F_{2}(1,1,2-t; 2,3; 1) = \frac{2(1 - \gamma - \psi(t+1))}{1-t}
\end{equation}
\noindent
to enable them demonstrate a new proof of (\ref{value-m0}). 

The graph of $M(a)$ shown in Figure 1, obtained by the numerical 
integration 
of (\ref{oloa-11}), has a well-defined {\em cusp} at $a = \ln 2$. In this 
paper, we provide analytic expressions for both branches of $M(a)$. The 
region 
$a > \ln 2$, determined in \cite{glasser-manna1}, is 
reviewed in this  section. The 
corresponding expressions for $0 < a  < \ln 2$ will be the 
content of the next section.

{{
\begin{figure}[ht]
\begin{center}
\includegraphics[width=3in]{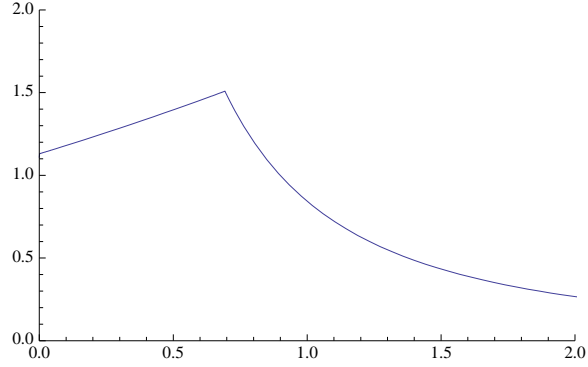}
\caption{The graph of $M(a)$ for $0 \leq a \leq 2$}
\label{figure-ma}
\end{center}
\end{figure}
}}

\section{The case $0 < a < \ln 2$} \label{sec-iden}
\setcounter{equation}{0}

Our starting point is the identity
\begin{equation}
M(a) = - \frac{e^{a}}{2 \pi} \imagpart{ 
\int_{0}^{1} e^{-at} \int_{-\pi}^{\pi} 
\frac{x (1+e^{ix})^{t}}{1-e^a + e^{ix} } \, dx \, dt},
\label{mofa-0}
\end{equation}
\noindent
as illustrated in \cite{glasser-manna1}. We outline the proof here
for sake of the reader's convenience. The identity 
\begin{equation}
\imagpart{ \frac{x}{ix + \ln \left[ 2 e^{-a} \cos x \right]} } = 
\frac{x^{2}}{x^{2} + \ln^{2} \left[ 2 e^{-a} \cos x \right]}
\end{equation}
\noindent
yields
\begin{equation}
M(a) = \frac{4}{\pi} 
\imagpart{ \int_{0}^{\pi/2} \frac{x \, dx }{ix + \ln \left[ 2 e^{-a} 
\cos x \right] } }. 
\end{equation}
\noindent
Under the assumption $a > \ln 2$, we have
\begin{equation}
\int_{0}^{\infty} e^{s \ln \left[ 2e^{-a} \cos x \right] + i x } \, ds = 
\frac{1}{ix + \ln \left[ 2e^{-a} \cos x \right]},
\end{equation}
\noindent
which implies
\begin{equation}
M(a) = \frac{2}{\pi} \imagpart{ \int_{-\pi/2}^{\pi/2} \int_{0}^{\infty} 
x e^{s( \ln \left[ 2 e^{-a} \cos x \right] + i x )} \, dx \, ds},
\end{equation}
\noindent
where one uses the fact that the imaginary part of the integrand is 
an {\em even} function of $x$. One more
identity
\begin{equation}
e^{isx} \cdot e^{s \ln \left[ 2 e^{-a} \cos x \right]} = 
e^{s \ln \left[ e^{-a} (1  + e^{2ix} ) \right]} 
\end{equation}
\noindent
and the change of variables $x \mapsto x/2$, give reason to 
\begin{equation}
M(a) = \frac{1}{2 \pi} \imagpart{ 
\int_{-\pi}^{\pi} \int_{0}^{\infty} x e^{s \ln \left[ e^{-a}(1+e^{ix}) \right]
}\, ds \, dx}.
\end{equation}
\noindent
Then evaluate the $s$-integral to obtain
\begin{equation}
M(a) = - \frac{1}{2 \pi} \imagpart{ \int_{-\pi}^{\pi} 
\frac{x \, dx }{\ln \left[ e^{-a} ( 1 + e^{ix} ) \right]}.
}
\label{mofa-18}
\end{equation}
\noindent
The formula
\begin{equation}
\frac{1}{\ln u} = \int_{0}^{1} \frac{u^{t} \, dt}{u-1}
\end{equation}
\noindent
now gives (\ref{mofa-0}) from (\ref{mofa-18}). 

\begin{note}
Even though the proof outlined here is valid for $a > \ln 2$, the identity
(\ref{mofa-0}) holds for $a > 0$. 
\end{note}

\noindent
{\bf Notation}: we use $b = e^{a}-1$ and assume $ 0 < a < \ln 2$, so that
$0 < b < 1$. \\

Expanding the terms $(1+e^{ix})^{t}$ and $1/(1-be^{-ix})$ 
in power series produces
\begin{equation}
M(a) = - \frac{e^{a}}{2 \pi} \int_{0}^{1} \int_{-\pi}^{\pi} xe^{-at} 
\sum_{j=0}^{\infty} \sum_{k=0}^{\infty} b^{j} \binom{t}{k} 
\sin [x(k-j-1)] \, dx \, dt.
\nonumber
\end{equation}
\noindent
The term corresponding to $k=j+1$ disappears and computing the 
$x$-integral we arrive at
\begin{eqnarray}
M(a) & = & e^{a} \int_{0}^{1} e^{-at} \sum_{j=0}^{\infty} 
\sum_{k=0}^{j} (-1)^{j-k} \frac{b^{j} \binom{t}{k}}{j+1-k} \, dt \label{star} \\
 & + & e^{a} \int_{0}^{1} e^{-at} \sum_{j=1}^{\infty} b^{j} \,
\sum_{\nu=1}^{j} \frac{(-1)^{\nu}}{\nu}  \binom{t}{\nu+j} 
\, dt. \nonumber 
\end{eqnarray}

\begin{lemma}
Let $t \in \mathbb{R}$ and $j \in \mathbb{N} \cup \{ 0 \}$. Then
\begin{equation}
\sum_{\nu=1}^{\infty} \frac{(-1)^{\nu}}{\nu} \binom{t}{\nu + j} = 
\binom{t}{j} \left[ \psi(j+1) - \psi(t+1) \right].
\end{equation}
\end{lemma}
\begin{proof}
The integral representation ($3.268.2$ in \cite{gr}):
\begin{equation}
\psi(p+1) - \psi(q+1) = - \int_{0}^{1} \frac{x^{p}-x^{q}}{1-x} \, dx,
\end{equation}
\noindent
yields
\begin{equation}
\psi(p+1) - \psi(q+1) = \sum_{j=1}^{\infty} (-1)^{j-1} 
\left( \binom{p}{j} - \binom{q}{j} \right). 
\end{equation}
Therefore the result is a consequence of the identity
\begin{equation}
\binom{t}{k}^{-1} \sum_{m=1}^{\infty} \frac{(-1)^{m}}{m} \binom{t}{m+k} - 
\sum_{m=1}^{\infty} \frac{(-1)^{m}}{m} \binom{t}{m}  = 
\sum_{m=1}^{k} \frac{(-1)^{m-1}}{m} \binom{k}{m}. 
\label{iden0}
\noindent
\end{equation}
\noindent
Apply the difference operator $\Delta a(k) = a(k+1)-a(k)$ and use
\begin{equation}
\binom{t}{k+1}^{-1} \binom{t}{m+k+1} - \binom{t}{k}^{-1} \binom{t}{m+k} =
-\frac{m}{k+1} \binom{t}{k+1}^{-1} \binom{t+1}{m+k+1}
\nonumber
\end{equation}
\noindent
to write the derived equation as 
\begin{equation}
- \frac{\binom{t}{k+1}^{-1}}{k+1} \sum_{m=1}^{\infty} (-1)^{m} \binom{t+1}
{m+k+1} = \Delta \sum_{m=1}^{k} \frac{(-1)^{m}}{m} \binom{k}{m}.
\label{iden-1}
\end{equation}
\noindent
The left hand side of (\ref{iden0}) reduces to $-1/(k+1)$ in view of the 
classical identity
\begin{equation}
\sum_{m=1}^{\infty} (-1)^{m-1} \binom{t+1}{m+k+1} = \binom{t}{k+1}.
\end{equation}
\noindent
A simple evaluation of the right hand side in 
(\ref{iden0}) also produces $-1/(k+1)$. 
We conclude that, up to a constant term with respect to the
index $k$, both sides of
(\ref{iden0}) are equal to the {\em harmonic number} $H_k$. 
The special case $k=0$ shows that this constant vanishes. 
\end{proof}

\medskip

Continuing from (\ref{star}), we thus have
\begin{eqnarray}
& & \label{mofa-88} \\ 
M(a) & = & e^{a} \int_{0}^{1} \sum_{j=0}^{\infty} b^{j} \, \sum_{k=0}^{j} 
\frac{(-1)^{j-k} \binom{t}{k} }{j+1-k} + 
\frac{e^{a}}{b} \int_{0}^{1} e^{-at} \sum_{j=1}^{\infty} b^{j} 
\binom{t}{j} \psi(j+1) \, dt \nonumber \\
& - & \frac{e^{a}}{b} \int_{0}^{1} e^{-at} \sum_{j=1}^{\infty} b^{j} 
\binom{t}{j} \psi(t+1) \, dt. \nonumber 
\end{eqnarray}

To simplify the first term in the previous expression observe
\begin{eqnarray}
\sum_{j=0}^{\infty} b^{j} \, \sum_{k=0}^{j} 
\frac{(-1)^{j-k} \binom{t}{k} }{j+1-k}  & = & 
 \sum_{k=0}^{\infty} b^{k} \binom{t}{k} \, \sum_{\nu=0}^{\infty} 
\frac{(-1)^{\nu} b^{\nu}}{\nu+1} \nonumber \\
& = & \frac{\ln(1+b)}{b} \sum_{k=0}^{\infty} \binom{t}{k}b^{k} 
 =  \frac{ae^{at}}{b}. \nonumber 
\end{eqnarray}
\noindent
We deduce that in (\ref{mofa-88}) the
first term  is  $a/(1-e^{-a})$. \\

The reduction of the second term in (\ref{mofa-88}) employs the following
result: 

\begin{lemma}
Let $ 0 < a < \ln 2$  and $t \in \mathbb{R}$. Then
\begin{equation}
\int_{0}^{1} e^{-at} \sum_{j=0}^{\infty} b^{j} \binom{t}{j} \psi(j+1) \, dt 
= \ln( 1 - e^{-a}) + \int_{1}^{\infty} \frac{e^{-at}}{t} \, dt.
\label{form-66}
\end{equation}
\end{lemma}
\begin{proof}
The Stirling numbers $s(j,k)$ satisfy
\begin{equation}
j! \binom{t}{k} = \sum_{k=0}^{j} | s(j,k) | t^{k},
\end{equation}
\noindent
so that
\begin{equation}
\int_{0}^{1} e^{-at} \sum_{j=0}^{\infty} b^{j} \binom{t}{j} \psi(j+1) \, dt 
= \frac{e^{-a} b \gamma}{a} + e^{-a} 
\sum_{j=1}^{\infty} (b^{j+1} \alpha_{j} - b^{j} \alpha_{j-1} ) \psi(j+1),
\end{equation}
\noindent
with
\begin{equation}
\alpha_{j}(a) = \frac{1}{j!} \sum_{k=0}^{j} \frac{|s(j,k)| k!}{a^{k+1}}.
\end{equation}
\noindent
The result now follows by summation by parts and the identity
\begin{equation}
\sum_{j=k}^{\infty} \frac{|s(j,k)| b^{j}}{j!} = 
\frac{\ln^{k}(1+b)}{k!}.
\end{equation}
\end{proof}

Therefore, the second term in (\ref{mofa-88}) is 
\begin{equation}
\mbox{ second term } = \frac{\ln(1-e^{-a})}{1-e^{-a}} + 
\frac{1}{1-e^{-a}} \int_{1}^{\infty} \frac{e^{-at}}{t} \, dt.
\end{equation}

\medskip

Finally, the third term in (\ref{mofa-88}) is
\begin{equation}
\mbox{third term}  =  - \frac{e^{a}}{b} \int_{0}^{1} e^{-at} 
\left( \sum_{j=1}^{\infty} \binom{t}{j} b^{j} \right) \psi(t+1) \, dt 
 =  - \frac{e^{a}}{b} \int_{0}^{1} ( 1  - e^{-at} ) \psi(t+1) \, dt. 
\nonumber 
\end{equation}

A direct computation shows that
$\int_{0}^{1} \psi(t+1) \, dt = 0$, and integration by parts gives
\begin{equation}
\mbox{third term } = \frac{a}{1-e^{-a}} \int_{0}^{1} e^{-at} \ln \Gamma(t+1) 
\, dt. 
\end{equation}
\noindent
The identity $\ln \Gamma(t+1) = \ln \Gamma(t) + \ln t$ now yields 
\begin{equation}
\mbox{third term } = \frac{a}{1-e^{-a}} 
\left( \int_{0}^{1} e^{-at} \ln t \, dt +  \int_{0}^{1} e^{-at} 
\ln \Gamma(t) \, dt \right). 
\nonumber 
\end{equation}
\noindent
Replacing (2.17), (2.19) and (2.21) into (\ref{mofa-88}) provides the  
following expression for $M(a)$: 
\begin{eqnarray}
M(a) & = & \frac{a}{1-e^{-a}} + \frac{\gamma}{a} + 
\frac{\ln(1 - e^{-a})}{1-e^{-a}} \label{int-77} \\
& + & \frac{a}{1-e^{-a}} \int_{0}^{\infty} e^{-at} \ln t \, dt + 
\frac{a}{1-e^{-a}} \int_{0}^{1} e^{-at} \ln \Gamma(t) \, dt.
\nonumber
\end{eqnarray}
\noindent
The term $\gamma/a$ comes from the index $j=0$ in the sum (\ref{form-66}). 
Next, we make use of (\ref{int-00}) to state our main result which, 
incidentally, is complementary to Theorem \ref{gm-1}. This 
settles a conjecture of 
O. Oloa stated in \cite{oloa1}.

\begin{theorem}
\label{mofa-small}
Assume $ 0 < a < \ln 2$. Then
\begin{eqnarray}
M(a)  =  \frac{\gamma}{a} + 
\frac{a + \ln(1-e^{-a}) - \gamma - \ln a}{1-e^{-a}}  +
\frac{a}{1-e^{-a}} \int_{0}^{1} e^{-at} \ln \Gamma(t) \, dt. 
\nonumber
\end{eqnarray}
\end{theorem}

This result enjoys a form similar to Theorem \ref{gm-1}.

\begin{corollary}
Assume $0 < a < \ln 2$. Then 
\begin{equation}
M(a) = \frac{\gamma}{a} + 
\frac{\left( a + \ln(1-e^{-a}) + \Gamma(0,a) \right)}{1-e^{-a}} +
\frac{1}{1-e^{-a}} \int_{0}^{1} e^{-at} \psi(t+1) \, dt.
\nonumber
\end{equation}
\noindent
where $\Gamma(0,a)$ is the incomplete gamma function.
\end{corollary}

\begin{proof}
Split up the first integral in (\ref{int-77}) and integrate by parts.
\end{proof}

Differentiating (\ref{mofa-small}) with respect to $a$ at $a=0$, and use
the classical value
\begin{equation}
\int_{0}^{1} \ln \Gamma(t) \, dt  = \frac{1}{2} \ln 2 \pi 
\end{equation}
\noindent
and also 
\begin{equation}
\int_{0}^{1} t \ln \Gamma(t) \, dt = 
\frac{\zeta'(2)}{2 \pi^{2}} + \frac{1}{6} \ln 2 \pi - \frac{\gamma}{12}
\end{equation}
\noindent
obtained in \cite{espmoll1}, produces
\begin{equation}
\int_{0}^{\pi/2} \frac{x^{2} \ln( 2 \cos x ) \, dx}{(x^{2} + \ln^{2}(2 \cos x) )^{2}} = \frac{7 \pi}{192} + \frac{\pi \, \ln 2 \pi}{96} - 
\frac{\zeta'(2)}{16 \pi}. 
\label{nice-1}
\end{equation}

Further differentiation of (\ref{mofa-small}) 
produces the evaluation of a family of integrals 
similar to (\ref{nice-1}). 

The integral in (\ref{mofa-small}) can be expressed in an alternative 
form. Define
\begin{equation}
\Lambda(z) := \lim\limits_{n \to \infty} \left( \sum_{j=1}^{n} 
\frac{j}{j^{2}+z^{2}}  - \ln n \right).
\end{equation}
\noindent
Observe that $\Lambda(0) = \gamma$, so $\Lambda(z)$ is  
a generalization of Euler's constant.

\begin{lemma}
Let $a > 0, \, c = 1 - e^{-a}$ and define $A = \ln 2 \pi + \gamma$. Then
\begin{equation}
\int_{0}^{1} e^{-at} \ln \Gamma(t) \, dt = \frac{A(a-c)}{a^{2}} - 
\frac{c}{2a} \Lambda \left( \frac{a}{2 \pi} \right) + 
2c \sum_{j=1}^{\infty} \frac{\ln j}{a^{2} + 4 \pi^{2} j^{2}}.
\end{equation}
\end{lemma}
\begin{proof}
Expand the exponential into a MacLaurin series and use the values
\begin{eqnarray}
\int_{0}^{1} t^{n} \ln \Gamma(t) \, dt & = & \frac{1}{n+1} 
\sum_{k=1}^{\lfloor{ \tfrac{n+1}{2} \rfloor}} (-1)^{k} \binom{n+1}{2k-1} 
\frac{(2k)!}{k(2 \pi)^{2k}} \left[ A \zeta(2k) - \zeta'(2k) \right] \nonumber \\
& - & \frac{1}{n+1} \sum_{k=1}^{\lfloor{ \tfrac{n}{2} \rfloor}}
(-1)^{k} \binom{n+1}{2k} \frac{(2k)!}{2(2 \pi)^{2k}} \zeta(2k+1) +
\frac{\ln \sqrt{2 \pi}}{n+1} \nonumber
\end{eqnarray}
\noindent
given  as $(6.14)$ in \cite{espmoll1}. Then calculate by interchanging
the resulting double sums. 
\end{proof}

The next corollary follows from the identity $M(a) = L(a) + \frac{\gamma}{a}$.

\begin{corollary}
Assume $0 < a < \ln 2$ and let $c= 1 - e^{-a}$. Then 
\begin{equation}
\int_{0}^{\infty} e^{-at} \ln \Gamma(t) \, dt = 
- \frac{\gamma + \ln a}{ace^{a}} + \frac{A(a-c)}{a^{2}c} -
\frac{1}{2a} \Lambda\left( \frac{a}{2 \pi} \right) + 
2 \sum_{j=1}^{\infty} \frac{\ln j}{a^{2} + 4 \pi^{2} j^{2}}. 
\end{equation}
\end{corollary}

\medskip

We now state two identities involving the function 
$f(t) = 2^{-t} \ln \Gamma(t)$. The proof of these identities was supplied to 
the authors by O. Espinosa. 

\begin{lemma}
\label{olivier}
The identities
\begin{eqnarray}
\int_{0}^{\infty} f(t) \, dt & = & 
2 \int_{0}^{1} f(t) \, dt - 
\frac{\gamma + \ln \ln 2}{\ln 2} \label{ident-1} \\
\int_{0}^{\infty} tf(t) \, dt & = & 2 \int_{0}^{1} (t+1) f(t) \, dt -
\frac{(\gamma + \ln \ln 2)(1 + 2 \ln 2) -1}{\ln^{2}2} \label{ident-2}
\end{eqnarray}
\noindent
hold.
\end{lemma}
\begin{proof}
The function $f(t)$ satisfies $f(t+1) = \frac{1}{2}f(t) + \frac{1}{2} 2^{-t}
\ln t$. Splitting the integral
\begin{equation}
\int_{0}^{\infty} f(t) \, dt = \int_{0}^{1} f(t) \, dt + 
\int_{0}^{\infty} f(t+1) \, dt
\end{equation}
\noindent
and using (\ref{int-00}) gives the first result. The proof of (\ref{ident-2}) 
is similar, it only requires differentiating (\ref{int-00}) with respect to
the parameter $a$.
\end{proof}

The reader will check that (\ref{ident-1}) is equivalent to the continuity of 
$M(a)$ at $a= \ln 2$. The identity (\ref{ident-2}) provides a proof of the 
last result is worthy of singular 
(pun-intended) interest.

\begin{theorem}
The jump of $M'(a)$
at $a = \ln 2$ is $4$.
\end{theorem}

\section{Conclusions} \label{sec-conclusions}
\setcounter{equation}{0}

The integral 
\begin{equation}
M(a):= \frac{4}{\pi} \int_{0}^{\pi/2} \frac{x^{2} \, dx }
{x^{2} + \ln^{2}( 2 e^{-a} \cos x ) },
\nonumber
\end{equation}
\noindent
satisfies
\begin{equation}
M(a)  = \frac{\gamma}{a} +
\int_{0}^{\infty} e^{-a t} \psi(t+1) \, dt,
\end{equation}
\noindent
for $a > \ln 2 $ and 
\begin{equation}
M(a) = \frac{\gamma}{a} + 
\frac{\left( a + \ln(1-e^{-a}) + \Gamma(0,a) \right)}{1-e^{-a}} +
\frac{1}{1-e^{-a}} \int_{0}^{1} e^{-at} \psi(t+1) \, dt
\nonumber
\end{equation}
\noindent
for $0 < a \leq \ln 2$.

\end{document}